\documentclass[12pt,twoside]{amsart}

\usepackage{enumerate}
\usepackage[all]{xy}
\usepackage[hyperindex=true,plainpages=false,colorlinks=false,pdfpagelabels]{hyperref}

\theoremstyle{plain}
\newtheorem{The}{Theorem}
\newtheorem*{The*}{Theorem}
\newtheorem{Pro}{Proposition}
\newtheorem{Lem}{Lemma}

\newtheorem*{Cor*}{Corollary}

\theoremstyle{definition}
\newtheorem*{Def}{Definition}

\theoremstyle{remark} 
\newtheorem{Rem}{Remark}
\newtheorem{Exa}{Example}
\newtheorem*{Rem*}{Remark}

\numberwithin{equation}{section}

\DeclareMathOperator{\End}{End}
\DeclareMathOperator{\Hom}{Hom}

\DeclareMathOperator{\SL}{SL}

\DeclareMathOperator{\su}{su}

\DeclareMathOperator{\Tr}{tr}             

\DeclareMathOperator{\vol}{vol}

\renewcommand{\Im}{\operatorname{Im}}

\newcommand{\dvector}[1]{{\left(\begin{matrix}#1\end{matrix}\right)}}
\DeclareMathOperator{\dbar}{\bar\partial}
\DeclareMathOperator{\del}{\partial}

\newcommand{\R}{\mathbb{R}}
\newcommand{\C}{\mathbb{C}}

\renewcommand{\H}{\mathbb{H}}

\begin{document}

\title[Higher genus minimal surfaces in $S^3$ and stable bundles]{Higher genus minimal surfaces in $S^3$ and stable bundles}

\author{Sebastian Heller}

\address{Sebastian Heller \\
  Institut f\"ur Mathematik\\
  Universit{\"a}t T\"ubingen\\ Auf der Morgenstelle
10\\ 72076 T¬ubingen\\ Germany
 }

\email{heller@mathematik.uni-tuebingen.de}


\date{\today}

\thanks{Author supported by SFB/Transregio 71}

\begin{abstract} 
We consider compact minimal surfaces $f\colon M\to S^3$ of genus $2$ which are homotopic to an embedding. We assume that the associated holomorphic bundle is stable. We prove that these surfaces 
can be constructed from a globally defined family of meromorphic  connections by the DPW method.
The poles of the meromorphic connections are at the Weierstrass points of the Riemann surface of order at most $2$. For the existence proof 
of the DPW potential we give a characterization of stable extensions $0\to S^{-1}\to V\to S\to 0$ of spin bundles $S$ by its dual $S^{-1}$ in terms of an associated element of $P H^0(M;K^2).$ We also consider the family of holomorphic structures associated to a minimal surface in $S^3.$ For surfaces of genus $g\geq2$ the holonomy of the connections is generically non-abelian and therefore the holomorphic structures are generically stable.
 \end{abstract}

\maketitle


\section{Introduction}
\label{sec:intro}
 The systematic investigation of harmonic maps from Riemann surfaces to $S^3$ (or more generally to symmetric spaces) has started with 
the introduction of the associated family of flat connections $\nabla^\zeta,$ 
see for example \cite{Po} or \cite{H1}. Using this family, there have been many deep results concerning harmonic $2-$spheres
and harmonic $2-$tori. For example, the space of CMC tori in $\R^3$ or minimal tori in $S^3$ is well-understood, see \cite{H1}, \cite{PS}, \cite{B} and \cite{KS}. 

On the other hand, there is no satisfactory treatment of compact
CMC or minimal surfaces of higher genus. Nevertheless, there is
a general method due to Dorfmeister, Pedit and Wu, \cite{DPW}, which produces, in principal, all such surfaces. 
The idea of the DPW method is to gauge $\nabla^\zeta$ into a family of meromorphic connections in a way which can be
reversed: It is shown that one can gauge the holomorphic family of connections by a positive gauge, i.e. a $\zeta-$depending family 
of endomorphisms of determinant $1,$
which is well-defined and upper-triangular with positive diagonal at $\zeta=0,$
into a family of meromorphic connections of a special form on simply connected domains. From such a family of meromorphic connections one obtains minimal surfaces as
follows: Take a $\zeta-$depending parallel frame and split it into the
unitary and the positive part by Iwasawa decomposition. Then the unitary part is a parallel frame of a family of unitary connections 
describing a minimal surface. The surface obtained in this fashion depends on the $\zeta-$depending starting condition of the parallel frame.
Dressing, i.e. changing this starting condition, will give new surfaces.

Making compact surfaces (or even surfaces with topology) via this 
method is much more complicated. One has to ensure, by choosing the family of meromorphic connections and the right dressing,
 that the frame of the unitary part is well-defined up to (unitary) holonomy, and that the surface closes.  This has been worked out only in very special cases, for example for trinoids, i.e. genus $0$ CMC surfaces with three Delauney ends, or tori. There are no examples for (closed) higher genus surfaces up to now.

The aim of this paper is to show that for compact oriented minimal surfaces in $S^3$ of genus $2$ it is possible to find a DPW potential with a nice behavior on the Riemann surface.
Our method applies to the following situations: We assume that the minimal surface is homotopic to an embedding.  This is equivalent to saying that the associated spin bundle $S$ has no holomorphic sections, see \cite{P}. The second assumption is that the holomorphic structure $(\nabla^\zeta)''$ at $\zeta=0$ is stable.
This condition is needed in technical details. For the only known example, the Lawson genus $2$ surface, it is satisfied. Moreover, the holomorphic structure $(\nabla^\zeta)''$ is stable for generic $\zeta\in\C$ as we show in section \ref{sec:family}. Under these assumptions the family of connections $\nabla^\zeta$ can be gauged globally to a family of meromorphic connections on $M$ with constant holomorphic structure given by the trivial extension of $S$ by $S^{-1}.$ The poles of the meromorphic connections are exactly at the Weierstrass points of the Riemann surface of order at most $2,$ see theorem \ref{main}.

In the first part of this paper, we recall the gauge theoretic description of minimal surfaces $f\colon M\to S^3$ in $S^3$ due to Hitchin \cite{H1}. We give explicit formulas of the occurring connections and sections in terms of geometric quantities like the Hopf field and the spinor connection.  We give a link to the local description of surfaces preferred by other authors.

In the second section, we introduce the associated family of flat connections
$\nabla^\zeta.$ We gauge this family by a $\zeta-$depending $B$
with special singularities such that the meromorphic connections have a constant holomorphic structure. This can be achieved by solving $\dbar$-equations on $M:$ The gauge $B$ must be a section in a bundle $E\to U\subset\C\ni\zeta$ whose fibers are finite dimensional spaces of holomorphic sections in a bundle over $M$ (varying in $\zeta$), see lemma \ref{solutionspace}.
The difficulty is in proving that the gauge has constant determinant $\det(B)=1$ in order to produce no more singularities (possibly varying in $\zeta$). The determinant is given by a map to a finite dimensional space $\det\colon E\to H^0(M;K^3),$ compare with lemma \ref{submersive}.
We use the implicit function theorem to find a gauge with $\det(B)=1.$ We reduce the proof of the surjectivity of the differential of $\det$ restricted to the $\zeta=0$ slice to some algebraic geometric condition
 to the holomorphic bundle given by $(\nabla^0)'',$ see lemma \ref{submersive} and theorem \ref{stability}.

To understand that algebraic geometric condition
we study non-trivial extensions $0\to S^{-1}\to V\to S\to 0$ over compact surfaces $M$ of genus $2$ such that $S$ has no holomorphic sections.
We define a natural $1:1$ correspondence between non-trivial extension of this form and elements of $P H^0(M;K^2).$ In this 
setup, we will identify the non-stable extensions with the bundles which do not satisfy the above mentioned condition (theorem \ref{stability}).

In the last chapter we consider compact oriented immersed minimal 
surfaces of genus $g\geq 2.$ We prove that the holonomy representation of $\nabla^\zeta$ is non-abelian and that the holomorphic bundle 
$(V,(\nabla^\zeta)'')$ associated to the connection $\nabla^\zeta$ is stable 
for generic $\zeta\in\C.$ This shows that the method of investigating the eigenline bundle of $\nabla^\zeta$ cannot work for $g\geq2$ as for tori. 

The author thanks Aaron Gerding, Franz Pedit and Nick Schmitt for helpful
discussions.

\section{Minimal Surfaces in $S^3$}
First we shortly describe a gauge-theoretic way of treating minimal surfaces in $S^3$ 
due to Hitchin \cite{H1}. 
We refer to \cite{LM} for details about Clifford 
algebras and Spinors, and to \cite{H1} for the main
source of the material described below. 

 Consider the round $3-$sphere $S^3$ with its tangent bundle trivialized
by left translation 
$$TS^3=S^3\times \Im\H$$ and Levi Civita connection given, with respect to the above trivialization, by $$\nabla=d+\frac{1}{2}\omega.$$
Here $\omega\in\Omega^1(S^3,\Im\H)$ is the Maurer-Cartan form of
$S^3$ which acts via adjoint representation on the Lie algebra $\Im\H=\su(2).$ This formula is
equivalent to the well-known characterization of the Levi-Civita connection by the property that for left-invariant vector 
fields $X,Y$ it satisfies $\nabla_XY=\frac{1}{2}[X,Y].$

  There are two equivalent 
complex representations of the spin group $S^3$ induced from the Clifford 
representation $$\C l(\R^3)=\H\otimes\C\oplus \H\otimes\C$$ on the
complex vector space $\H$ with complex structure given by right 
multiplication with $i.$  It is well-known that
$S^3$ has an unique spin structure.
We consider the associated complex
spin bundle $$V=S^3\times\H$$
with complex structure given by right multiplication with $i\in\H.$
We have a complex hermitian metric $(.,.)$ on it given by the trivialization and by the identification
 $\H=\C^2.$
The Clifford multiplication is given by
$$TS^3\times V\to V; (\lambda, v)\mapsto \lambda v$$
where $\lambda\in\Im \H$ and $v\in\H.$ This is clearly complex linear.
The unitary spin connection is given by
\begin{equation}\label{nablaS3}
\nabla=\nabla^{spin}=d+\frac{1}{2}\omega,
\end{equation}
where the $\Im\H-$valued Maurer-Cartan form acts by left multiplication on the quaternions. Via this construction the tangent bundle
$TS^3$ identifies as the skew hermitian 
trace-free complex linear endomorphisms of $V.$
 
Let $M$ be a Riemann surface and $f\colon M\to S^3$ be a conformal 
immersion. Then the pullback $\phi=f^*\omega$ of the Maurer-Cartan form satisfies the structural equations
$$d\phi+\frac{1}{2}[\phi\wedge\phi]=0.$$
Another way to write this equation is
\begin{equation}\label{structure}
d^\nabla\phi=0,
\end{equation}
where $\nabla=f^*\nabla=d+\frac{1}{2}\phi,$ with $\phi\in\Omega^1(M;\Im\H)$ acting via adjoint representation. Conversely, every solution 
$\eta\in\Omega^1(M;\Im\H)$ to $d^{\nabla^\eta}\eta=0,$
where $\nabla^\eta=d+\frac{1}{2}\eta,$ gives rise to
a map $f\colon \hat M\to S^3$ from the universal covering $\hat M$ of $M$ unique up to translations in $S^3.$ 

From now on we only consider the case of $f$ being minimal. Under the assumption of $f$ being
conformal $f$ is minimal if and only if it is harmonic. This is
exactly the case when
\begin{equation}\label{harmonic}
d^\nabla*\phi=0.
\end{equation}
Consider $\phi\in\Omega^1(M;f^*TS^3)\subset
\Omega^1(M;\End_0(V))$ via the interpretation of
$TS^3$ as the bundle of trace-free skew hermitian
endomorphisms of $V.$ Decompose
$\frac{1}{2}\phi=\Phi+\bar \Psi$
into $K$ and $\bar K$ parts, i.e. 
$\Phi=\frac{1}{2}(\phi-i*\phi)\in\Gamma(K\End_0(V))$ 
and $\bar\Psi=\frac{1}{2}(\phi+i*\phi)\in\Gamma(\bar K\End_0(V)).$
The property of $\phi$ being skew symmetric translates to
$\bar\Psi=-\Phi^*,$ i.e.
$$\frac{1}{2}\phi=\Phi-\Phi^*.$$
Then $f$ is conformal if and only if $\Tr\Phi^2=0,$ see \cite{H1}.
In view of a rank $2$ bundle $V$ and $\Tr\Phi=0$ this is equivalent to
\begin{equation}\label{detPhi}
\det\Phi=0.
\end{equation}
Note that $f$ is an immersion if and only if $\Phi$ is nowhere vanishing. In other words, the branch points of $f$ are exactly the 
zeros of $\Phi.$
Moreover the equations
\ref{structure} and \ref{harmonic} are equivalent to
\begin{equation}\label{dbarPhi}
\dbar^\nabla\Phi=0,
\end{equation}
where
$\dbar^\nabla=\frac{1}{2}(d^\nabla+i*d^\nabla)$ 
is the induced holomorphic structure on $KV\to M.$

Of course equation \ref{dbarPhi} does not contain
the property $\nabla=d+\frac{1}{2}\phi,$ 
i.e. that $\nabla-\frac{1}{2}\phi$ is trivial on $V.$ 
Locally, or on simply 
connected sets, this is equivalent to the flatness of
$\nabla-\frac{1}{2}\phi,$ which is the same as the following formula
\begin{equation}\label{curvature}
F^\nabla=[\Phi\wedge\Phi^*],
\end{equation}
as one easily computes. 

Conversely, given an unitary rank $2$ bundle $V\to M$ over
a simply connected Riemann surface with special unitary connection
$\nabla$ and trace free field $\Phi\in\Gamma(K\End_0(V))$ without 
zeros, which satisfy the equations \ref{detPhi}, \ref{dbarPhi} and \ref{curvature}, we 
get a conformally immersed  minimal surface as follows:
By equation \ref{dbarPhi} and \ref{curvature}, the unitary connections
$\nabla^L=\nabla-\Phi+\Phi^*$ and $\nabla^R=\nabla+\Phi-\Phi^*$
are flat. Because $M$ is simply connected they are gauge equivalent.
Due to the fact that $\Tr\Phi=0,$ the determinant bundle $\Lambda^2V$ is trivial with respect to all these connections. Hence, the gauge 
is $SU(2)=S^3-$valued with differential $\phi=2\Phi-2\Phi^*.$ Thus it is 
a conformal immersion. The harmonicity follows from equation \ref{dbarPhi}.

\subsection{The Spinor Bundle of a Minimal Surface}\label{Spinor_Bundle}
We describe how the spin structure of the immersion
$f\colon M\to S^3$ can be seen in this setup. The geometric significance
of the spin structure is described in Pinkall
\cite{P}, see also the literature therein. We give formulas
which relate the data on $V$ obtained in the
previous part to the data usually used 
to describe a surface, for example the Gauss map and
 the Hopf field. Again, this part is based on \cite{H1}. 
 
In this section we consider the bundle $V$ with its
holomorphic structure $\dbar:=\nabla''.$
As we have seen the complex part $\Phi$ of the differential of a 
conformal minimal surface satisfies $\Tr\Phi=0$ and $\det\Phi=0,$
but is nowhere vanishing.  We obtain a well-defined
holomorphic line subbundle
$$L:=\ker\Phi\subset V.$$
Because $\Phi$ is nilpotent the image of $\Phi$ 
satisfies $\Im\Phi\subset K\otimes L.$ Consider 
the holomorphic section
$$\Phi\in H^0(M;\Hom(V/L,KL))$$
without zeros. The holomorphic structure $\dbar-\Phi^*$ turns $V\to M$ 
into the holomorphically trivial bundle $\underline\C^2\to M.$ As $\Tr\Phi^*=0,$ the determinant line bundle $\Lambda^2V$ of $(V,\dbar)$ is holomorphically trivial. This implies $V/L=L^{-1}$ and we obtain
$$\Hom(V/L,KL)=L^2K$$ as holomorphic line bundles. Because
$L^2K$ has a holomorphic section $\Phi$ without zeros, we get
$$L^2=K^{-1}.$$
Hence, its dual bundle $S=L^{-1}$ is a spinor bundle of the 
Riemann surface $M.$ Clearly, $S^{-1}$ is the only $\Phi-$invariant
line subbundle of $V.$ 

There is another way to obtain the bundles $S$ and $S^{-1}$ which provides a link to the quaternionic 
holomorphic geometry, see \cite{BFLPP}.
 Let $R\colon M\to\Im\H$ be the normal of $f$
with respect to the trivialization of the tangent bundle 
$TS^3=S^3\times\Im\H.$ Here, $R$ stands for 
the right normal vector when considering the surface as lying in
$S^3\subset\H.$ As we have seen we 
can consider $V$ as the trivial quaternionic line bundle 
$\underline\H\to M.$ Note that scalar multiplications with quaternions is 
from the right in order to commute with the Clifford multiplication.
We define a complex quaternionic linear structure $\mathcal J$ by
$$v\mapsto -Rv.$$ 
Note that $\mathcal J$ can be seen as the operator given by
Clifford multiplication with the negative of the determinant, i.e.
for a positive oriented orthonormal basis $X,Y\in TM$ we have
$\mathcal J(v)=Y\cdot X\cdot v,$ where $\cdot$ is Clifford 
multiplication. Then $V$ splits into $\pm i$ eigenspace of $\mathcal J:$
\begin{equation}\label{decompositionV}
V=E\oplus\bar E:=\{v\in V\mid \mathcal J=v i\}\oplus\{v\in V\mid \mathcal J=-v i\}.
\end{equation}
This decomposition is orthogonal with respect to the (complex) unitary
metric  on $V=\underline \H.$
\begin{Pro}
The kernel $S^{-1}$ of $\Phi$ is given by the $-i$ eigenspace of $\mathcal J.$
\end{Pro}
\begin{proof}
It is sufficient to prove $\bar E\subset\ker\Phi.$ Note that for all 
$v\in\underline\H$ the vector $v+\mathcal Jv=v-Rvi$ is an element of
$\bar E.$ With $4\Phi=\phi-i*\phi,$ the proof is simply a matter of 
computation. 
\end{proof}
We have seen that there exists a holomorphic subbundle $S^{-1}$ of 
$(V,\dbar),$ and that the determinant line bundle $\Lambda^2V$ is trivial. Therefore, $(V,\dbar)$ is a nontrivial extension of
$S$ by $S^{-1}:$
$$0\to S^{-1}\to V\to S\to0.$$
This means that with respect to the decomposition 
$V=S^{-1}\oplus S$ the holomorphic structure $\dbar$ can be written as
\begin{equation}\label{dbaronV}
\dbar=\dvector{&\dbar^{spin^*} & \bar b \\ &0& \dbar^{spin}}
\end{equation}
where 
$\bar b\in\Gamma(\bar K\Hom(S,S^{-1})=\Gamma(\bar KK^{-1}),$
and $\dbar^{spin}$ and $\dbar^{spin^*}$ are the holomorphic 
structures on $S$ and $S^{-1},$ respectively. It is well-known (\cite{H1}) that there exists a relation between $\bar b$ and the Hopf differential $Q$
of the minimal surface. We want to determine the exact form of this relation. More generally, we want to find out geometric formulas for the connection
$\nabla$ on the pullback $V\to M$ of the spinor bundle of $S^3.$

To do so recall that the pullback of the Levi-Civita connection $\nabla$
splits, with respect to the decomposition $f^*TS^3=TM\oplus\R$
into tangential and normal part, into
$$\nabla=\dvector{& \nabla^M  & -II^* \\  &II & d }$$
with $\nabla^M$ being the Levi-Civita connection on $M.$ Here $II$ is the second fundamental form of the surface which is a symmetric bilinear form $II\in\Gamma(T^*M\otimes T^*M).$ The Weingarten operator is given by $A=-II^*\in\End(TM),$ i.e.
$<A(X),Y>=-II(X,Y)$ for tangent vectors $X,Y\in TM.$ 
The $K^2-$part of the second fundamental form is called 
the Hopf field $Q.$
For minimal surfaces in $S^3$ it is holomorphic, i.e. $Q\in H^0(M;K^2).$ Its zeros are exactly the umbilics of the surface. 

The connection
$$\tilde\nabla=\dvector{& \nabla^M  & 0 \\  &0 & d }$$
is a $SO(3)-$connection, too, and it induces a 
unitary connection
$\tilde\nabla$ on the spinor bundle $V.$ But it reduces to a $SO(2)-$connection,
the Levi-Civita connection on $M,$ so the corresponding connection
on $V=S^{-1}\oplus S$ is the spin connection of $M.$ Its $\bar K-$part is given by
$$\tilde\nabla''=\dvector{&\dbar^{spin^*} & 0\\ & 0 & \dbar^{spin}}.$$
The difference of these two connections on $V$ is 
the trace-free skew adjoint operator
\begin{equation}\label{nablas}
\nabla-\tilde\nabla=\frac{1}{2}(\mathcal J\circ A)\cdot.
\end{equation}
This means for all $X\in TM,$ $\psi\in\Gamma(V)$ we have
$\nabla_X\psi-\tilde\nabla_X\psi=\frac{1}{2}(\mathcal J A(X))\cdot\psi$
with $\cdot$ being the Clifford multiplication and
$A$ the Weingarten operator. Note that this difference
is an off-diagonal endomorphism. One can compute that its $K-$part
vanishes on $S,$ and as an operator in
$K\Hom(S^{-1},S)=K^2$ it is exactly $-\frac{i}{2}Q.$
The adjoint $Q^*$ of $Q\in H^0(K^2)$ with respect to the hermitian product  $(.,.)$ is determined by $$(Q(X)v,w)=(v,Q^*(X)w)$$
for all $X\in TM,$ $v\in S^{-1},$ and $w\in S.$ This gives a well-defined 
section $Q^*\in \Gamma(\bar KK^{-1}).$
As $\frac{1}{2}\mathcal J\circ A$ is skew adjoint, the 
extension class $[\bar b]\in H^1(M;K^{-1})$ of $V$ is given by the representative
\begin{equation}\label{extensionclassV}
\bar b=-\frac{i}{2}Q^*\in\Gamma(\bar KK^{-1}).
\end{equation}
Because $Q$ is holomorphic, one can deduce that the extension class
$[\bar b]\in H^1(M,K^{-1})$ is non-zero (or $Q=0,$ which corresponds
to a totally geodesic $2-$sphere), see \cite{H1}. Altogether we obtain
\begin{Pro}\label{connection_data}
Let $f\colon M\to S^3$ be a conformal minimal immersion with associated complex unitary rank $2$ bundle $(V,\nabla).$ Let $V= S^{-1}\oplus S$ be the unitary decomposition, where $S^{-1}=\ker\Phi\subset V$ and $\Phi$ is the $K-$part of the differential of $f.$
With respect to this decomposition the connection can be written as
$$\nabla=\dvector{&\nabla^{spin^*} & -\frac{i}{2} Q^*\\ & -\frac{i}{2} Q & \nabla^{spin}},$$
where $\nabla^{spin}$ is the spin connection corresponding to the Levi-Civita connection on $M$ and $Q$ is the Hopf field of $f.$

The Higgsfield $\Phi\in H^0(M,K\End_0(V))$ can be identified with $$\Phi=1\in H^0(M;K\Hom(S,S^{-1})),$$ and its adjoint $\Phi^*$ is given by the volume form $\vol$ of the induced Riemannian metric.
\end{Pro}

\subsection{Local description}\label{local}
Next we give a link of the gauge theoretic description of minimal surfaces in $S^3$ with the local treatment of CMC surfaces in $\R^3$ or $S^3.$ The later is usually used by people working with the DPW
method. Moreover, the construction of minimal surfaces out of a meromorphic potential uses the local description: The Iwasawa decomposition of a (local) parallel frame of the meromorphic connection,
which is after a trivialization given by the meromorphic potential, splits out a so called extended frame $\mathcal F$ depending on $\zeta.$ This extended frame is nothing but the frame of the family of connections $\nabla^\zeta$ with respect to a corresponding  trivialization.

Let $U\subset M$ be a simply connected open subset and $z\colon U\to\C$ be a holomorphic chart. Write $g=e^{2u}|dz|^2$ for a function $u\colon U\to\R.$
Choose a local holomorphic section $s\in H^0(U;S)$ with $s^2=dz,$ and let $t\in H^0(U,S^{-1})$ be its dual holomorphic section. 
Then $$(e^{-u/2}t,\ e^{u/2} s)$$ is a special unitary frame of $V=S^{-1}\oplus S$ over $U.$ Write the Hopf field $Q=q (dz)^2$ for some local holomorphic function $q\colon U\to\C.$ 

The Levi-Civita connections of conformally equivalent metrics
$g=e^{2\lambda}g_0$ and $g_0$ differ on the canonical bundle $K$ by
the the form $-2\del\lambda=-(d\lambda-i*d\lambda)\in\Gamma(K).$ Therefore, the connection form of $\nabla^{spin}$ with respect to the local frame $s$ is given by
$-\del u,$ and with respect to $e^{u/2} s,$ it is given by $\frac{1}{2}i*du.$ From proposition \ref{connection_data} the connection form of $\nabla$ with respect to 
$( e^{-u/2}t,\ e^{u/2} s)$ is
$$\dvector{&-\frac{1}{2}i*du & -\frac{i}{2} e^{-u}\bar q d\bar z\\ & -\frac{i}{2} e^{-u}q d z& \frac{1}{2}i*du}.$$
The Higgsfield $\Phi$ and its adjoint $\Phi^*$ are given by 
$$\Phi=\dvector{&0 &  e^{u}d z\\ & 0& 0},\,\,\,\,\,\, \Phi^*=\dvector{&0 & 0\\ &  e^{u}d \bar z& 0}$$
with respect to the frame $( e^{-u/2}t,\ e^{u/2} s).$
These formulas are well-known, see \cite{DH}, or, in slightly other notation, \cite{B}. Therefore, the associated family of flat connections, see equation \ref{family}, takes locally 
the same form which is used to compute minimal surfaces in $S^3$ out of a meromorphic potential.

\section{DPW: From minimal surfaces to meromorphic connections}
We restrict our considerations to compact oriented minimal surfaces in
$S^3$ of genus $2.$ There is some hope that surfaces of genus $g\geq 3$ can be treated similar to surfaces of genus two, but there will be
much more technical difficulties in general.
We assume that the surface is homotopic to an embedding, and that the holomorphic bundle $(V,\nabla'')$ is stable.
We use the notations of the previous section.

 From equations \ref{curvature} and \ref{dbarPhi} we see that the curvature 
 of 
 \begin{equation}\label{family}
 \begin{split}
 \nabla^\zeta:=\nabla+\zeta^{-1}\Phi-\zeta\Phi^*
 \end{split}
 \end{equation}
 vanishes for all $\zeta\in\C\setminus\{0\}.$ The connections are special unitary for $\zeta\in S^1\subset\C,$ and $\SL(2,\C)-$connections for $\zeta\in\C.$
This family of connections plays a very important role in the theory of harmonic maps, as one might see from \cite{Po}, \cite{H1}, \cite{B}, \cite{DPW}, and \cite{KS}, or others.

The DPW method, see \cite{DPW}, shows that, on simply connected domains, every family of connections $\nabla^\zeta$ of the form \ref{family} can be obtained from meromorphic data, namely the DPW potential. We will not describe the details of this construction, one might consult \cite{DH} or \cite{DPW}. 
The idea is, that a local parallel frame $\Psi$ (with respect to a trivialization) of the meromorphic connection can be split by Iwasawa decomposition
$$\Psi=\mathcal F B$$
into a unitary part $\mathcal F$ and a positive part $B.$ The positive part 
will be the singular gauge described below (in the trivialization). The unitary part $\mathcal F$ is the parallel frame (in a corresponding trivialization) of a family of connections $\nabla^\zeta$ of the form \ref{family} obtained for minimal surfaces.

We will prove here that for compact oriented minimal surface in $S^3$ of genus $2,$ under the conditions described above, one can gauge the holomorphic family of connections $\nabla^\zeta$ globally to a family of meromorphic connections on $M.$ The $\dbar-$part of this meromorphic family is constant and given by the trivial extension of $S$ by $S^{-1}.$
We prove that the poles of the connections are exactly at the Weierstrass points of the Riemann surface of order at most $2.$
\begin{Rem*}
The condition that the surface is homotopic to an embedding translates to the property that the spin bundle $S\to M$ has no
holomorphic sections, see \cite{P}.
\end{Rem*}

\begin{Def}
A meromorphic connection on a holomorphic vector bundle $(V,\dbar)$ over a Riemann surface is a connection with 
singularities which can be written with respect to a local holomorphic frame as $d+\xi,$ where $\xi$ is an meromorphic endomorphism-valued $1-$form.
\end{Def}

Of course, meromorphic connections are flat on surfaces. For line bundles $L\to M$ there is a class of meromorphic connections which are in $1:1$ correspondence with meromorphic sections of $L$ by declaring the section to be parallel. Moreover there exists the degree formula
$$res(\nabla)=-deg(L)$$ on Riemann surfaces, where $res(\nabla)$ is the sum of all local (well-defined)
residua (of the locally defined connection forms).

The condition that $\nabla^\zeta\cdot B,$ the gauge of $\nabla^\zeta$ by $B,$ is a holomorphic family of meromorphic connections
on the holomorphic bundle $S^{-1}\oplus S$ translates easily to
\begin{equation}\label{dbarB}
\dbar^{spin}B=(\frac{i}{2}Q^*+\zeta\Phi^*) B,
\end{equation}
where $Q^*\in\Gamma(\bar KK^{-1})$ and $\Phi^*\in\Gamma(\bar KK)$ are given as in the previous section, and $\dbar^{spin}$ is the holomorphic structure of the endomorphism bundle of the direct sum bundle $S^{-1}\oplus S.$

With respect to the unitary splitting $V=S^{-1}\oplus S$ we
write $$B=b_\zeta=\dvector{ a & b\\ c & d}$$ 
with $a,d\in\Gamma(M\times U;\underline\C),\, b\in\Gamma(M\times U;K^{-1})$ and
$c\in\Gamma(M\times U;K),$ where $U\subset\C$ is a small neighborhood of $0$ in the $\zeta-$plane.
Then equations \ref{dbarB} are equivalent to
\begin{equation}\label{gauging}
\begin{split}
\dbar a&=\frac{i}{2} Q^*c\\
\dbar c&=\zeta\Phi^*a\\
\dbar b&=\frac{i}{2} Q^*d\\
\dbar d&=\zeta\Phi^*b,\\
\end{split}
\end{equation}
where the $\dbar-$operators are the obvious ones on the trivial holomorphic bundle $\underline\C,$ on the canonical bundle $K$ and on its dual $K^{-1}.$
One cannot solve these equations globally (on compact surfaces) without singularities.
For example, if one makes a $\zeta-$expansion of $a,b,c,d$ and starts with $d_0=1,$ then, by Serre duality, there does not exist a solution to
$\dbar b_0=\frac{i}{2} Q^*1,$ because $[Q^*]\in H^1(M,K^{-1})$ is non-zero. What one has to do is to allow singularities of the following kind: 
\begin{Def}
Let $W\to M$ be a holomorphic vector bundle over a Riemann surface.
 A local section $s\in\Gamma(U\setminus\{p\},W)$ has a pole-like singularity of order $k$ at $p\in U$ if it can locally be written as $\frac{t}{z^k}$ for a locally non-vanishing section $t$ and a holomorphic chart $z$ centered at $p.$
 The space of all sections with pole-like singularities will be denoted by
$\hat\Gamma(M,W).$
\end{Def}
To solve the $\dbar$ problem at hand, we allow pole-like singularities
at the Weierstrass points of the Riemann surface.
More concrete, we 
take two divisors $D\neq\tilde D$ with 
$$L(D)=L(\tilde D)=K S.$$
In fact one can take $D=Q_1+Q_2+Q_3$ and $\tilde D=Q_4+Q_5+Q_6,$ where $Q_1,..,Q_6$ are the Weierstrass points of the Riemann surface (in the right order corresponding to the spin
structure $S$). To see this, take
a Weierstrass point $Q_1.$ Then there exists two uniquely determined points $P_1,P_2\in M$ such that
$L(P_1+P_2-Q_1)=S$ by Riemann-Roch. Since $S$ has no holomorphic sections, we see $L(P_1+P_2)\neq K,$ and one easily 
obtains that $P_1=Q_2$ and $P_2=Q_3$ are Weierstrass points.
It is  clear that $Q_1,$ $Q_2$ and $Q_3$ must be pairwise disjoint.
Now take another Weierstrass point $Q_4$ and the corresponding Weierstrass points $Q_5,Q_6$ such that $L(Q_5+Q_6-Q_4)=S.$ Again, it is clear that
$Q_1,..,Q_6$ are pairwise disjoint.

Now, we multiply with $s_D$ and $s_{\tilde D}$ to guarantee 
the existence of solutions: We consider sections
$$\tilde a=a\otimes s_D\in\Gamma(M\times U;KS),\, \tilde c=c\otimes s_D\in\Gamma(M\times U;K^2S),$$ and $$\tilde b=b\otimes s_{\tilde D}\in\Gamma(M\times U;S),\, \tilde d=d\otimes s_{\tilde D}\in\Gamma(M\times U;KS).$$
In order to solve the equations \ref{gauging} $\tilde a\oplus\tilde c$ and $\tilde b\oplus\tilde d$ have to be holomorphic sections of the holomorphic bundles described in the following
\begin{Lem}\label{solutionspace}
There exists $\epsilon>0$ and a holomorphic bundle
$V_1\to B(0;\epsilon)\subset\C$ of rank $6$ with the property
$V_1(\zeta):=\ker(\dbar_1^\zeta),$
where
$$\dbar_1^\zeta=\dvector{\dbar & -\frac{i}{2}Q^* \\ -\zeta \Phi^* & \dbar} \colon\Gamma(M;KS\oplus K^2S)\to\Gamma(M;\bar KKS\oplus \bar KK^2S).$$
Similarly, there exists a holomorphic bundle $V_2\to B(0;\epsilon)\subset\C$ of rank $2$ with the property
$V_2(\zeta):=\ker(\dbar_2^\zeta),$
where
$$\dbar_2^\zeta=\dvector{\dbar & -\frac{i}{2}Q^* \\ -\zeta \Phi^* & \dbar} \colon\Gamma(M;S\oplus KS)\to\Gamma(M;\bar KS\oplus \bar KKS).$$
\end{Lem}
\begin{proof}
Note that for a (holomorphic) family of elliptic operators the minimal kernel dimension is attained on an open set. Over this open set the kernel bundle is holomorphic. For details see \cite{BGV} or \cite{BPP}.

It remains to prove that $\ker(\dbar_1^0)$ and $\ker(\dbar_2^0)$
have minimal dimension. By Riemann-Roch
$$index(\dbar_1^\zeta)=6.$$ 
With Serre duality one obtains that
$$\ker(\dbar_1^0)\cong H^0(M;KS)\oplus H^0(M;K^2S)$$
has dimension $6.$ 

Similarly,
$$index(\dbar_2^\zeta)=2,$$ and
$$\ker(\dbar_2^0)\cong H^0(M;KS)$$
has dimension $2.$
\end{proof}
From lemma \ref{solutionspace} we see that we can find a gauge
$B$ as follows: Take a holomorphic section
$$\tilde B=(B_1,B_2)\in H^0(B(0;\epsilon),V_1\oplus V_2).$$
Then $$B:=(B_1\otimes s_{-D}, B_2\otimes s_{-\tilde D})$$ is a 
($\zeta$-depending) section of $\End(V)$ with pole-like singularities at $Q_1,..,Q_6.$ If we can
choose $\tilde B$ such that $B$ has constant determinant $\det(B)=1,$ then
$B$ gauges $\nabla^\zeta$ into a holomorphic family of meromorphic connections with constant $\dbar-$part given by the trivial extension of $S$ by $S^{-1}.$

In order to ensure $\det(B)=1,$ we have to study the following determinant:
\begin{Lem}\label{submersive}
The determinant map $$\det\colon V_1(0)\oplus V_2(0)\to H^0(M,K^3)$$ given by $\det(\dvector{a\\ c},\dvector{ b \\ d}):= ad-bc$
has surjective differential at the point
$\dvector{s_D & s_{\tilde D}'\\ 0 & s_{\tilde D}},$ where $s_{\tilde D}'\in\Gamma(M;S)$
is the unique solution of $\dbar s_{\tilde D}'=\frac{i}{2}Q^* s_{\tilde D},$ exactly in the case that $(V,\nabla'')$ is stable.
\end{Lem}
\begin{proof}
First of all $\det$ maps to $H^0(M;K^3)$ as a consequence of the
equations \ref{gauging}.
So it is well-defined and holomorphic. 

Its
derivative $d_p\det$ at $p:=\dvector{s_D & s_{\tilde D}'\\ 0 & s_{\tilde D}}$ is given by the map
\begin{equation}\label{differential}
\begin{split}
\dvector{\alpha+q' & \beta' \\ q & \beta}\mapsto
s_D\beta+ s_{\tilde D} \alpha+s_{\tilde D} q'-q s_{\tilde D}'
\end{split}
\end{equation}
for $\alpha,\beta\in H^0(M;KS),$ $q\in H^0(M;K^2S)$ and
solutions $q'\in\Gamma(M;KS)$ of $\dbar q'=\frac{i}{2}Q^* q$ and
$\beta'\in\Gamma(M;S)$ of $\dbar\beta'=\frac{i}{2}Q^*\beta.$
Note that $\dim H^0(M;KS)=2,$ $\dim H^0(M;K^2S)=4,$ and
$\dim H^0(M;K^3)=5.$ 

Looking at the zeros one sees that
$s_D\beta+ s_{\tilde D} \alpha=0$ exactly in the case that
$\beta=\lambda  s_{\tilde D}$ and $\alpha=-\lambda s_D$ for some
$\lambda\in\C.$ Therefore, the differential $d_p\det$ maps the subspace given by  $\{q=0,\ q'=0\}$ to a $3-$dimensional subspace of $H^0(M;K^3).$

Consider a basis $(q_1,..,q_4)$ of $H^0(M;K^2S)$
with divisors given by
\begin{equation}
\begin{split}
(q_1)&=D+2Q_1\\
(q_2)&=D+2Q_4\\
(q_3)&=\tilde D+2Q_1\\
(q_4)&=\tilde D+2Q_4.\\
\end{split}
\end{equation}
It can be easily seen that the differential $d_p\det$ maps $q_3$ and $q_4$
into the $3-$dimensional subspace described above. Any element in the subspace spanned by $q_1,q_2$ can be written as $\omega s_D$
for some $\omega\in H^0(M;K).$ Its image lies in the $3-$dimensional subspace of $H^0(M;K^3)$ exactly in the case that
there exists $\alpha,\beta\in H^0(M;KS)$ such that
\begin{equation}\label{deco}
\omega (s_Ds_{\tilde D}'-s_D's_{\tilde D})=s_{ D}\beta+s_{\tilde D}\alpha.
\end{equation}
The decomposition \ref{deco} is possible for nonzero $\omega$ exactly in the case of a non-stable
bundle $(V,\nabla''),$ see Theorem \ref{stability} and Remark \ref{Deco1}.
\end{proof}
\begin{The}\label{main}
Let $\nabla^\zeta$ be the holomorphic family of flat connections (\ref{family}) on $V$ associated to a compact oriented immersed minimal surface $f\colon M\to S^3$ of genus $2.$ Assume that $(V,\nabla''=(\nabla^0)'')$ is stable and that $f$ is homotopic to an embedding. Let $S$ be the associated spinor bundle of $f.$ Then, there exists an order
$Q_1,..,Q_6$ of the six Weierstrass points of $M$ 
such that $KS=L(Q_1+Q_2+Q_3)= L(Q_4+Q_5+Q_6).$
 
There exists holomorphically $\zeta-$dependent gauge
$$B:\zeta\in \tilde B(0;\epsilon)\subset\C\to \hat\Gamma (\End(V))$$ with pole-like singularities at $Q_1,..,Q_3$ up to order $1$ 
in the first column (with respect to the unitary decomposition
$V=S^{-1}\oplus S)$ and
pole like singularities at $Q_4,..,Q_6$ up to order $1$ 
in the second column
such that
$\det B_\zeta=1$ for all $\zeta$ and such that the gauged connection
$$\hat\nabla^\zeta:=\nabla^\zeta\cdot B_\zeta$$ is a holomorphic family of meromorphic connections $\hat\nabla^\zeta$ for 
$\zeta\in B(0;\epsilon)\setminus\{0\}\subset\C$ on the (fixed) direct sum 
holomorphic vector bundle $S^{-1}\oplus S.$ Moreover, $B(0)=\dvector{1 & *\\ 0 & 1},$ which implies that $B$ is a positive gauge.

The connections $\hat\nabla^\zeta$ have poles up to order $1$ on the diagonal 
(with respect to the unitary decomposition 
$V=S^{-1}\oplus S)$
at $Q_1,..,Q_6$ and poles up to order $2$ in the lower left entry
at $Q_1,..,Q_3$ and in the upper right at $Q_4,..,Q_6.$
The family $\hat\nabla^\zeta$ has an expansion in $\zeta$ of the form
$$\hat\nabla^\zeta=\dvector{\nabla_{0}^* & \zeta^{-1}+\omega \\ -\frac{i}{2}Q & \nabla_{0}}+higher\ order\ terms,$$
where $\nabla_0$ is a meromorphic connection on $S,$ $\omega\in\mathcal M(M;\C),$ and $Q\in H^0(K^2)$ is the Hopf field of the minimal surface.
\end{The}
\begin{proof}
Because of Lemma \ref{submersive} and the implicit function theorem we can find locally around 
$\zeta=0$ a holomorphic section $(B_1,B_2)$
of $V_1\oplus V_2\to B(0;\epsilon)$ with $B_1(0)=s_{Q_1+Q_2+Q_3}$ and
$\det(B_1,B_2)(\zeta)=s_{Q_1+Q_2+Q_3}s_{Q_2+Q_4+Q_6}$
for small $\zeta.$ Here $\det$ is defined on $V_1\oplus V_2$ analog as in Lemma \ref{submersive}. Then the gauge given by
$$B:=(B_1\otimes s_{-Q_1-Q_2-Q_3},B_2\otimes s_{-Q_4-Q_5-Q_6})$$ is of the desired form.

The expansion of the family of connections $\hat\nabla^\zeta$
and its pole behavior can be easily computed. 
\end{proof}
By this theorem one knows what kind of DPW potential one should use to construct genus $2$ minimal surfaces $f\colon M\to S^3.$ With respect to a meromorphic trivialization of $S^{-1}\oplus S$ the DPW potential
$\xi$ takes values (for each $\zeta\in\C^*$) in an explicitly known finite dimensional vector space (depending only on $M$). Of course, this theorem does not give new informations about the behavior of the
potential $\xi$ into the $\zeta-$direction. 
\begin{Rem}
The condition on the stability, being essential in the proof presented here, is natural in some sense. First of all, the only known example in genus $2,$ Lawson's genus $2$ surface (\cite{L}), has a stable holomorphic bundle $(V,\nabla''),$ see \cite{He} or example \ref{exa} below. 
Moreover, 
stability is an open condition, compare with theorem \ref{stability} and section \ref{sec:family}.
\end{Rem}
\begin{Rem}
In principle it should be possible to prove a theorem of this kind for all compact surfaces of higher genus. One of the main problems
will be to find out such detailed informations about the poles
of the meromorphic connections.
\end{Rem}
\begin{Rem}
The theorem \ref{main} also applies for compact oriented CMC surfaces in $\R^3$ or $S^3$ of genus $2,$ as one sees from the discussion in section \ref{local}.
\end{Rem}
\begin{Exa}
In concrete situations one should be able to find out more 
informations about the meromorphic connections, for example
in the case that the surface has many symmetries. This has 
already been done by the author (\cite{He}) in the case of Lawson's genus $2$ surface: The Riemann surface is given by the equation
$$y^3=z^4-1.$$
The Weierstrass points are the points $Q_1,Q_2,Q_3$ lying over $0$
and $Q_4,Q_5,Q_6$ lying over $\infty,$ and the spinor bundle is $S=L(Q_1+Q_2-Q_3).$ The umbilics $P_1,..,P_4$ of the minimal surface are the branch points.
In a meromorphic trivialization of $S^{-1}\oplus S$ 
given by the meromorphic sections 
 $s=s_{Q_4+Q_5+Q_6-P_1-P_2-P_3-P_4}\in\mathcal M(M,S^{-1})$ and 
 $t=s_{-Q_6-Q_5-Q_6+P_1+P_2+P_3+P_4}\in\mathcal M(M,S),$
the connection form of $\hat\nabla^\zeta$ is given by
$$\dvector{-\frac{4}{3}\frac{z^3}{z^4-1}+\frac{A}{z}  & \zeta^{-1}+Bz^2 \\ \frac{G}{(z^4-1)} +\frac{\zeta H}{z^2(z^4-1)}& \frac{4}{3}\frac{z^3}{z^4-1}-\frac{A}{z}}dz.$$
Here $A,B,G,H$ are $\zeta-$depending holomorphic functions well-defined at $\zeta=0$ which satisfy $H= A+A^2$ and $B= -\frac{1}{G}(-\frac{1}{3}+A+(\frac{1}{3}-A)^2).$
\end{Exa}

\section{Stable extensions $0\to S^{-1}\to V\to S\to 0$}\label{sec:stable}
For general informations and details about extensions
and stable bundles we refer to \cite{NR}. 
We restrict to the case that $S$ is a spinor bundle over a compact Riemann surface of genus $2$ which has no holomorphic sections.

It is well-known that non-trivial extensions 
$$0\to S^{-1}\to V\to S\to 0$$ correspond to elements of
$PH^1(M;K^{-1})$ as follows: any two sections $b_1,b_2\in\Gamma(M;\bar KK^{-1})$ give rise to holomorphic isomorphic extensions
via $\dbar=\dvector{\dbar & b_k\\0 & \dbar}$ on $S^{-1}\oplus S$ 
if and only they are in the same class in $PH^1(M;K^{-1}).$
\begin{The}\label{classifying}
There exists a projective isomorphism
$$\Phi\colon PH^1(M;K^{-1})\to PH^0(M;K^2),$$
which does only depend on the spin bundle $S.$ Hence, the space
$PH^0(M;K^2)$ classifies nontrivial extensions $0\to S^{-1}\to V\to S\to 0.$
\end{The}
\begin{proof}
Consider a section $b\in\Gamma(M;\bar KK^{-1})$ which defines
a non-trivial extension $0\to S^{-1}\to V\to S\to 0.$
Let $\alpha,\beta\in H^0(M;KS)$ be a basis of $H^0(M;KS)$ and $\alpha',\beta'\in\Gamma(M;S)$ be the unique solutions of $\dbar\alpha'=b\alpha$
and $\dbar\beta'=b\beta.$ Then
$$\mathcal Q:=\alpha'\beta-\alpha\beta'\in H^0(M;K^2).$$
Clearly, the line $\C\mathcal Q\in PH^0(M;K^2)$ does not depend on the chosen basis $\alpha,\beta.$ Moreover, if we consider 
$\tilde b=b+\dbar X$ for $X\in\Gamma(M;K^{-1})$ we see that
$\tilde\alpha'=\alpha'+X\alpha$ and $\tilde\beta'=\beta'+X\beta«$
are the corresponding solutions. Hence $\C\mathcal Q$ does only depend on the class $[b]\in H^1(M;K^{-1}).$

Because the solutions clearly depend linearly on $b\in\Gamma(M;\bar KK^{-1})$ it remains to show that $\mathcal Q\neq 0$ for $0\neq[b]\in H^1(M;K^{-1}).$ To see this note that $\alpha$ and $\beta$ have no
common zeros. Therefore, if $\mathcal Q$ would vanish, the solution $\beta'$
would have zeros at the zeros of $\beta.$ 
Here, and later on in this section, we say that a smooth section
$s$ has a zero at $p$ of order $k$ if and only if $s\otimes (s_{-p})^k$ is smooth.
Hence, there 
would be a solution of $$\dbar t=b$$ for $$t (=\beta'\otimes \beta^{-1})\in\Gamma(M;K^{-1})=\Gamma(M;S\otimes (KS)^{-1}).$$ Here,
$\beta^{-1}\in\mathcal M(M;(KS)^{-1})$ is the dual meromorphic section of $\beta.$ By Serre duality and the non-vanishing of
$[b]\in H^1(M;K^{-1})$ there cannot be a solution $t.$
\end{proof}
\begin{Rem*}
One should not mistake the line $\C\mathcal Q$ in the space of holomorphic quadratic differentials associated to a non-trivial extension and the Hopf field $Q$ of a minimal surface. But in the case of the Lawson genus $2$ surface, the Hopf differential generates the line associated to the holomorphic structure $\nabla'',$ see example \ref{exa}. 
\end{Rem*}
The advantage of the description given by theorem \ref{classifying} is the following
\begin{The}\label{stability}
A nontrivial extension $0\to S^{-1}\to V\to S\to 0$ over a compact Riemann surface of genus $2$ with corresponding $\C\mathcal Q\in P H^0(M;K^2),$ such that $S$ has no holomorphic sections, 
is stable if and only if there
exist $0\neq\omega\in H^0(M;K)$ and $\alpha,\beta\in H^0(M;KS)$ 
such that $$\omega\mathcal Q=\alpha\beta.$$
\end{The}
\begin{Rem}\label{Deco1}
Note that each element $\psi$ of the $3-$dimensional subspace 
$W\subset H^0(M;K^3)$ spanned by products of two holomorphic sections in $H^0(M;KS)$ can be written as a product
$\psi=\alpha\beta$ for holomorphic sections 
$\alpha,\beta\in H^0(M;KS).$
\end{Rem}
\begin{proof}
A nontrivial extension $0\to S^{-1}\to V\to S\to 0$ given by
$b\in\Gamma(M;\bar KK^{-1})$ is non-stable if and only if
there exists a point $P\in M$ such that
$$b^\perp=\{q\in H^0(M,K^2)\mid q(P)=0\},$$
where $$b^\perp=\{q\in H^0(M,K^2)\mid \int_M (b,q)=0\},$$
see lemma 5.2 of \cite{NR}.

We need to characterize the zeros of $\mathcal Q:=\Phi([b]).$
Note that each holomorphic quadratic differential is the product
of two holomorphic differentials.
Let $(\mathcal Q)=P_1+..+P_4$ for $P_k\in M,$ such that 
$L(P_1+P_2)=L(P_3+P_4)=K.$ For each point $P\in M$ there exists an unique pair of points $\tilde P,\hat P\in M$ such that
$KS=L(P+\tilde P+\hat P).$ We claim that $P$ is a zero of $\mathcal Q$
if and only if 
$$\int_M(b,q)=0$$
for all $q\in H^0(M;K^2)$ with $q(\tilde P)=q(\hat P)=0,$ counted with multiplicities (only important for the case of $\tilde P=\hat P)$. 
Because
$L(\hat P+\tilde P)\neq K$ the space of $q\in H^0(M;K^2)$ with $q(\tilde P)=q(\hat P)=0$ is $1-$dimensional, and it is determined by $P.$ For $k=1,..,4$ we consider a basis $(s=s_{P_k+\tilde P_k+\hat P_k},\ t)$ of 
$H^0(M;KS).$ Because $s$ and $t$ have no common zeros,
the solution $s'$ of $\dbar s'=b s$ has a zero at 
$P_k,$ too. Therefore, there exists a solution 
$$\varphi(=s'\otimes s_{-P_k})\in\Gamma(M;SL(-P_k))$$ of
$$\dbar \varphi=b s_{\tilde P_k+\hat P_k}=b s\otimes s_{-P_k}\in\Gamma(M;\bar KSL(-P_k)).$$
But by Serre duality, there exists such a solution if and only if
$$\int_M(bs_{\tilde P_k+\hat P_k},\omega)=0$$ for all
$\omega\in H^0(M;SL(P_k)).$ This space is $1-$dimensional
because $S$ has no holomorphic sections.
Then  $s_{\tilde P_k+\hat P_k}\omega\in H^0(M; KS L(-P_k)SL(P_k))$ is a holomorphic 
quadratic differential which spans the $1-$dimensional space
$$\{q\in H^0(M;K^2)\mid q(\tilde P_k)=q(\hat P_k)=0\}.$$ This proves the assertion for the zeros of $\mathcal Q.$  

There exists a holomorphic differential $\omega$ and
$\alpha,\beta\in H^0(M;KS)$ with $\omega\mathcal Q=\alpha\beta$ 
if and only if $\{\tilde P_1,\hat P_1\}\cap\{P_3,P_4\}$ is non-empty.
This can be easily deduced from the facts that $L(\tilde P_1+\hat P_1)\neq K$ and that $S$ has no holomorphic sections. 

If $\{\tilde P_1,\hat P_1\}\cap\{P_3,P_4\}$  is non-empty, we can assume that $\tilde P_1=P_3.$ Then $\{\tilde P_3,\hat P_3\}
=\{P_1,\hat P_1\}.$ Let us first assume that $P_1\neq P_3.$
By the characterization of the zeros of $\mathcal Q$ this implies that
$$\int_M(b,q)=0$$ for all
$q\in H^0(M;K^2)$ with $q(\tilde P_1)=q(\hat P_1)=0$ or $q(P_1)=q(\hat P_1)=0.$ But these holomorphic quadratic 
differentials span the $2-$dimensional space 
$$\{q\in H^0(M,K^2)\mid q(\hat P_1)=0\},$$ and we see that the 
extension is non-stable. If $P_1=P_3,$ then $P_2=P_4,$ too, and we
have $L(2 P_1+\hat P_1)=KS.$ With the same methods as above
one can show that the property that $\mathcal Q$ has a zero of order $2$ at $P_1$ implies that $\int_M(b,q)=0$ for all holomorphic quadratic differentials with $q(\hat P_1)=0.$ Again, this implies that the extension is non-stable.

Conversely, assume that the extension is non-stable. Therefore,
there exists a point $P\in M$ such that $\int_M(b,q)=0$ for all
holomorphic quadratic differentials with $q(P)=0.$ By the characterization of the zeros of $\mathcal Q$ one easily sees that
$\tilde P$ and $\hat P$ are zeros of $\mathcal Q.$ First assume that
$\tilde P\neq\hat P.$ Let $\omega$ be a non-zero holomorphic 
differential with $\omega(P)=0.$ Then $$D:=(\omega)+(\mathcal Q)-P-\tilde P-\hat P$$ is a positive divisor with $L(D)=KS,$ and 
$\omega\mathcal Q=s_{P+\tilde P+\hat P} s_D$ is a decomposition as 
required. Now assume $\tilde P=\hat P.$ Because $\int_M(b,q)=0$ for all holomorphic differentials with $q(P)=0,$ one can show that $\tilde P$ is a zero of order $2$ of $\mathcal Q.$ Let $\omega$ be a non-zero holomorphic 
differential with $\omega(P)=0.$ Again
$$D:=(\omega)+(\mathcal Q)-P-2\tilde P$$ is a positive divisor with 
$L(D)=KS,$ and $\omega\mathcal Q=s_{P+2\tilde P} s_D$ is a 
decomposition as required.
\end{proof}
\begin{Exa}\label{exa}
We claim that the line $\C\mathcal Q\in PH^0(M,K^2)$ associated to the holomorphic structure of Lawson's genus $2$ surface, see\cite{L}, is given by its Hopf differential 
$Q.$ Let $P_1,..,P_4$ be the umbilics of the surface, i.e. the zeros of $Q.$ They correspond to the points lying over $0$ and $\infty$ in the hyper-elliptic picture 
$$y^2=z^6-1$$
of the Riemann surface. Then $\omega_1=\frac{1}{\sqrt{z^6-1}}dz,$ $\omega_2=\frac{z}{\sqrt{z^6-1}}dz,$ is a basis of the space of holomorphic differentials. As in \cite{He}, 
$Q$ is given by a multiple of $\omega_1\omega_2,$ and $Q^*$ is perpendicular to $(\omega_1)^2$ and $(\omega_2)^2$ as a consequence of the symmetries of the Lawson surface. By the proof of theorem \ref{stability} the zeros of $\C\mathcal Q$ are the zeros of the Hopf differential. Moreover, one sees from \cite{NR} or from the characterization of theorem \ref{stability}, that the holomorphic structure $\nabla''$ is stable.
\end{Exa}

\section{The family of holomorphic structures}\label{sec:family}
In the last section we consider immersed compact oriented minimal surfaces in $S^3$ of genus $g\geq 2.$
We will prove that the holomorphic structure
$$\dbar^\zeta:=(\nabla^\zeta)''=\dvector{\dbar^{spin*} & -\frac{i}{2}Q^*\\ \zeta\Phi^* & \dbar^{spin}}$$
on $V$ is stable for generic $\zeta\in\C.$ We need
\begin{Pro}
Any holomorphic subbundle $L$ of degree $0$ of $(V,\dbar^\zeta)$
for $\zeta\in S^1\subset\C^*$ is parallel with respect to
$\nabla^\zeta.$
\end{Pro}
\begin{proof}
The holomorphic bundle $L^*$ has a unique unitary flat connection $\nabla^{L^*}.$ Then, the subbundle $L\subset V$ gives rise to a holomorphic section $i\in H^0(M;L^*\otimes V).$ 
We denote the induced flat unitary connection on $L^*\otimes V$
by $\nabla=\del+\dbar.$
As in 
\cite{H1} we obtain from flatness
$$\dbar\del i=0$$
and
$$\int_M(\del i,\del i)=\int_M(\dbar\del i,i)=0.$$ Thus $i$ is parallel, 
and $L$ is a parallel subbundle of $(V,\nabla^\zeta).$
\end{proof}
If $(V,\dbar^\zeta)$ would not be stable for generic $\zeta\in\C$ the
holonomy of $\nabla^\zeta$ would be abelian for all $\zeta\in\C:$ For $\zeta\in S^1\subset\C$ this follows easily from the fact that with $L\subset V$ parallel, also $Lj\subset V$ is parallel. Because 
the holonomy depends holomorphically on $\zeta,$ this implies the assertion. But the following theorem shows that this is not possible for $g\geq 2.$
\begin{The}
The holonomy representation of $\nabla^\zeta$ is non-abelian for generic $\zeta\in\C^*.$ As a consequence, $(V,\dbar^\zeta)$ is stable for generic $\zeta\in\C^*.$
\end{The}
\begin{proof}
The proof is based on the
observation that many arguments of Hitchin (\cite{H1}) remain true 
for higher genus surfaces under the assumption of 
abelian holonomy: 
First of all one would obtain a splitting 
$$V=L_\zeta\oplus L^*_\zeta$$
into parallel subbundles of $\nabla^\zeta$ for $\zeta$ in a punctured neighborhood $\hat U$ of $\zeta=0.$ Of course, this decomposition is holomorphic in $\zeta$ locally, but it might be that the subbundles
interchange as $\zeta$ goes around $0.$ We can also assume that $L_\zeta^2$ is not holomorphically trivial for $\zeta\in\hat U.$ As in \cite{H1} we get a holomorphic decomposition of the trace free endomorphisms
$$\End_0(V,\dbar^\zeta)=L_\zeta^2\oplus \underline\C\oplus L_\zeta^{-2},$$
where the $\underline\C$-part corresponds to (trace free) diagonal endomorphisms corresponding to the decomposition $V=L_\zeta\oplus L^*_\zeta.$ From this one sees $\dim H^0(M,\End_0(V,\dbar^\zeta))=1$
for $\zeta\in\hat U.$ Moreover, a generator of this $1-$dimensional space is parallel with respect to $\nabla^\zeta.$ The bundle
$$H^0(M,\End_0(V,\dbar^\zeta))\to\hat U$$
extends to $\zeta=0.$ Consider a local trivializing section $\Psi,$ i.e. a holomorphic family of holomorphic
sections $$\zeta\in\C\mapsto \Psi_\zeta\in H^0(M,\End_0(V,\dbar^\zeta))$$
which is non-vanishing for small $\zeta.$ This section is covariant
constant with respect to
$\nabla^\zeta$ for small $\zeta\neq0.$
Expanding $\Psi_\zeta=\Psi^0+\zeta\Psi^1+...$ around
 $\zeta=0$ implies
 $$[\Psi^0,\Phi]=0.$$ 
 This yields that $\Psi^0$ is a (non-zero) 
 holomorphic section in 
 $\Hom(S,S^{-1})\subset \End_0(V,\dbar^\nabla).$
 But $\deg\Hom(S,S^{-1})=2-2g<0$ for $g\geq2,$ 
 and we obtain a contradiction.
\end{proof}
Because of this theorem it is not possible to define an eigenline spectral curve which does not depend on the chosen generator $\gamma\in\{\alpha_1,\beta_1,..\alpha_g,\beta_g\}\subset\pi^1(M).$ The eigenlines for different $\gamma$ do not coincide, and the whole machinery which was so successful for tori cannot be applied for higher genus $g\geq2.$ 


\begin{thebibliography}{10}

\bibitem[BGV]{BGV} Berline, N., Getzler, E., Vergne, M., {Heat Kernels and Dirac Operators}, Grundlehren, Springer-Verlag, 2003.

\bibitem[B]{B} Bobenko, A. I., {\em Surfaces of constant mean curvature and integrable equations}, translation in Russian Math. Surveys 46, no. 4, 1991.

\bibitem[BFLPP]{BFLPP} Burstall, F. E., Ferus, D., Leschke, K.
  Pedit, F., and Pinkall, U., {\em Conformal geometry of surfaces in
    $S\sp 4$ and quaternions}, Lecture Notes in Mathematics 1772,
  Springer-Verlag, Berlin, 2002.

\bibitem[BLPP]{BLPP}  Bohle, C., Leschke, K., Pedit, F., and Pinkall, U.,
{\em Conformal maps from a 2-torus to the 4-sphere},
arXiv:0712.2311v1 [math.DG].

\bibitem[BPP]{BPP}  Bohle, C., Pedit, F., and Pinkall, U.,
{\em The spectral curve of a quaternionic holomorphic line bundle over a $2-$torus},
arXiv:0904.2475v1 [math.DG].

\bibitem[DH]{DH} Dorfmeister, J., Haak, G., {\em Meromorphic potentials and smooth surfaces of constant mean curvature},  Math. Z.  224,  no. 4, 1997.

\bibitem[DPW]{DPW} Dorfmeister, J., Pedit, F., Wu, H.,  {\em Weierstrass type representation of harmonic maps into symmetric spaces},  Comm. Anal. Geom.  6,  no. 4, 1998.

\bibitem[FLPP]{FLPP} Ferus, D., Leschke, K., Pedit, F., Pinkall, U, {\em Quaternionic holomorphic geometry: Pl\"ucker formula, Dirac eigenvalue estimates and energy estimates of harmonic $2$-tori}, Invent. Math. 146, no. 3, 2001.


\bibitem[GHPS]{GHPS} Gerding, A., Heller, S. Pedit, F., and Schmitt, N., {\em Global aspects of integrable surface geometry}, Proceedings for "Integrable Systems and Quantum Field Theory at Peyresq, Fifth Meeting.

\bibitem[GriHa]{GriHa} Griffith, P., and Harris, J., {\em Principles of
  algebraic geometry}, Pure and applied mathematics, John Willey \& Sons, New
  York, 1978.
  
  
  \bibitem[He]{He} Heller, S. {\em Lawson's genus $2$ surface and DPW}, in preperation.
 
 \bibitem[H1]{H1} Hitchin, N. J., {\em Harmonic maps from a $2$-torus to the $3$-sphere}, J. Differential Geom. 31, no. 3, 1990. 
 
  \bibitem[H2]{H2} Hitchin, N. J., {\em The self-duality equations on a Riemann surface}, Proc. London Math. Soc. (3) 55, no. 1, 1987. 
  
\bibitem[KS]{KS} Kilian, M., Schmidt, M. U.,{\em On the moduli of constant mean curvature cylinders of finite type in the $3-$sphere},
arXiv:0712.0108v2 [math.DG].

\bibitem[KPS]{KPS} Karcher, H., Pinkall, U., Sterling, I., {\em  New minimal surfaces in $S\sp 3$}, J. Differential Geom. 28 , no. 2, 1988.


\bibitem[L]{L} Lawson, H. B.,{\em Complete minimal surfaces in $S\sp{3}$}, Ann. of Math. (2), 92, 1970 .

\bibitem[LM]{LM} Lawson, H.B., and Michelsohn, M.L., 
{\em Spin Geometry}, Princeton, 1990.


\bibitem[NS]{NS} Narasimhan, M. S., Seshadri, C.S.,{\em Stable and unitaty vector bundles on a compact Riemann surface}, Ann. of Math. (2) 82, 1965.

\bibitem[NR]{NR} Narasimhan, M. S., Ramanan, S.,{\em Moduli of vector bundles on a compact Riemann surface}, Ann. of Math. (2) 89, 1969.





\bibitem[P]{P} Pinkall, U., {\em Regular homotopy classes of immersed surfaces}, Topology 24 , no. 4, 1985.   

\bibitem[PS]{PS} Pinkall, U., Sterling, I., {\em On the classification of constant mean curvature tori}, Ann. of Math. (2), 130, no. 2, 1989.

\bibitem[Po]{Po} Pohlmeyer, K., {\em Integrable Hamiltonian systems and interactions through quadratic constraints},  Comm. Math. Phys.  46, no. 3, 1976.



\end{thebibliography}
\end{document}